\documentclass[11pt]{amsart}
\usepackage{amsfonts,amsthm,amsfonts,amssymb,amstext,amsopn,amscd}
\usepackage{enumerate,color}
\usepackage{showidx}
\usepackage{latexsym,enumerate}
\usepackage{amsmath,amsthm,amsopn,amstext,amscd,
amsfonts,amssymb,mathrsfs}
\usepackage{graphicx}
\usepackage{tikz}
\usepackage[english]{babel}
\usepackage[latin1]{inputenc}
\usepackage{times}
\usepackage[T1]{fontenc}
\usepackage{color}
\usepackage{times}
\usepackage{verbatim}
\usetikzlibrary{arrows,shapes}
\numberwithin{equation}{section}
\usepackage[labelformat=empty]{caption}
\def\Box{\vrule height 6pt depth 1pt width 4pt}
\date{}
\begin{document}
\begin{center}
{\Large Symmetry and stability of non-negative }
\\[0.1cm]
{\Large solutions to degenerate elliptic equations in a ball }
\\[0.5cm]
{\sc  F. Brock} $^1$ 
{\sc \& P. Tak\'{a}\v{c}} $^1 $
 
\setcounter{footnote}{1}
\footnotetext{University of Rostock, Institute of Mathematics, 18057 Rostock, Ulmenstr. 69, Haus 3, Germany,  e-mail:  friedemann.brock@uni-rostock.de, peter.takac@uni-rostock.de  }
\end{center}
\small
\vspace*{0.5cm}
{\bf Abstract: }
We consider non-negative distributional solutions $u\in C^1 (\overline{B_R } )$ to the equation \\
$-\mbox{div} [g(|\nabla u|)|\nabla u|^{-1} \nabla u ]  = f(|x|,u)$ in a ball $B_R$, with $u=0$ on $\partial B_R $, where $f$ is continuous and non-increasing in the first variable and $g\in C^1 (0,+\infty )\cap C[0, +\infty )$, with $g(0)=0$ and $g'(t)>0$ for $t>0$. According to a result of the first author, the solutions satisfy a certain 'local' type of symmetry.
Using this, we first prove that the solutions are radially symmetric  
provided that $f$ satisfies appropriate growth conditions near its zeros. 
\\ 
In a second part we study the autonomous case, $f=f(u)$.
The solutions of the equation are critical points for an associated variation problem. We show under rather mild conditions that global and local minimizers of the variational problem are radial.
%and second variations of $J$ in certain directions 
\\[0.3cm]
\\[0.3cm]
\\[0.3cm]
\normalsize
\small
\renewcommand{\baselinestretch}{1.1}
\normalsize
\section*{1. Introduction}
\renewcommand{\theequation}{\thesection.\arabic{equation}}
\setcounter{section}{1} 
\setcounter{equation}{0}
Gidas, Ni and Nirenberg have proved in their celebrated paper 
\cite{GNN} that positive solutions to some uniformly elliptic problems in symmetric domains are  symmetric. A sample result is the following
\\[0.1cm]
{\bf Theorem A. } 
{\sl Let $B_R $ be a ball in $\mathbb{R}^N $ with radius $R>0$   
 centered at the origin. Further, let  
$f\in C[0, +\infty )$ and 
$$
f= f_1 + f_2 ,
$$
where $f_1 $ is Lipschitz continuous and $f_2 $ is non-decreasing.  
Finally,  
let $u$ be a solution of the following problem
$$
{\bf (P_0 )}
\left\{ 
\begin{array}{l}  
u\in C^2  (\overline{B_R }) ,\\ 
 u> 0, \quad
 -\Delta  u =f( u) 
\ \mbox{ in } \ B_R ,
\\
 u=0 \ \mbox{ on } \ \partial B_R .
\end{array}
\right.
$$ 
Then $u$ is radially symmetric and radially decreasing. More precisely, there is a function $U\in C^1 ([0, +\infty )) $ such that 
\begin{equation}
\label{radU}
  u(x) =U(|x|)  \quad \forall 
x\in B_R , \quad \
\ \mbox{and }\ 
 U'(r) <0 \quad \forall r\in (0, R] .
\end{equation} 
}
The proof of this result was based on  the  {\sl Alexandrov-Serrin moving plane method } 
which turned 
out to be a very powerful technique in showing symmetry results for solutions to elliptic and parabolic problems in the sequel, see \cite{Se}, \cite{GNNrn}, \cite{BeNi}, 
\cite{LiNi1}, \cite{LiNi2}, 
\cite{DamPaRa}, 
\cite{SeZ} and \cite{DamRa}. The survey article \cite{weimingni} and the monograph \cite{Fra} provide many further references about this subject. 
\\
Let us point out that extensions of 
Theorem A  
to the case that  $f= f(|x|,u)$ and $f$ is non-increasing in the first variable have also been treated in \cite{GNN}. 
\\
Further, one can show radial symmetry  of solutions if one replaces $\Delta $   by the $p$-Laplace operator when $1<p<2$, see \cite{DamPa1}, \cite{DamPa2}. 
\\
However, the following example shows that the question of symmetry becomes more delicate if one replaces the Laplacian by  a degenerate elliptic operator, and/or if the assumptions for $f$ are weakened.
\\[0.1cm]
{\bf Example 1.1.}  (see also \cite{Br2})
\\
Let $p,s \in \mathbb{R}$, with $p>1 $ and $s> \frac{p}{p-1}$, and define  
\begin{eqnarray*}
w(x) := 
\left\{ 
\begin{array}{ll}
(1-|x|^2 )^s & \mbox{ if $ |x|\leq 1 $}
\\
0 & \mbox{ if $|x|>1$}
\end{array}
\right.
,
\\
v(x) := \left\{ 
\begin{array}{ll}
1 & \mbox{ if $ |x|< 5 $}
\\
1- \left(
\frac{|x|^2 -25}{11} 
\right)^s
  & \mbox{ if $5\leq |x|\leq 6$}
\end{array}
\right.
.
\end{eqnarray*}
Then we choose $x^1 , x^2 \in B_4 $ with $|x^1 -x^2 |\geq 2 $ and set
$$
u(x):= v(x)+ w(x-x^1 )+  w(x-x^2 ) \quad \forall x\in B_6 .
$$
The graph of $u$ is built by three radially symmetric mountains, one of them having a plateau at height $1$ while the other two are congruent to each other with their feet lying on the plateau. 
\\[0.1cm]
\includegraphics[scale=1.0]{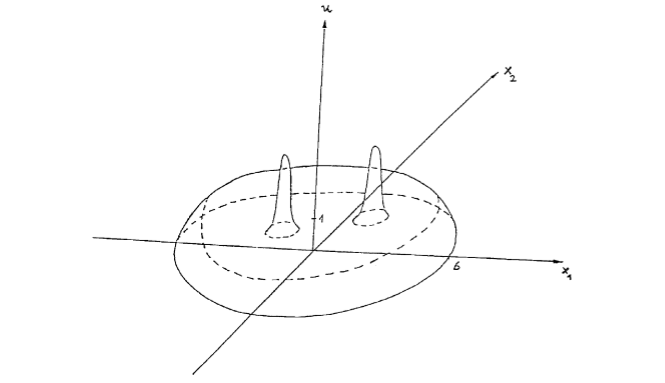}
\\[0.1cm]
After a short computation we see that $u$ is a weak solution of
the following problem
$$
{\bf (P_1 )} 
\left\{
\begin{array}{l}
 u\in C^1 (\overline{B_6 }) ,
\\ 
u>0, \quad -\Delta _p u \equiv \, -\mbox{div} \left( |\nabla u|^{p-2} \nabla u \right) = f(u) \ \mbox{in }\ B_6 ,
\\
 u=0 \ \mbox{on }\ \partial B_6 ,
\end{array}
\right.
$$
where
\begin{equation*}
f(u) := \left\{ 
\begin{array}{l}
\left( 
\frac{2s}{11} 
\right) 
^{p-1} 
\left[ 
25 +11 
\left( 
1-u
\right) 
^{1/s } 
\right] 
^{(p/2)-1} 
\cdot 
\left( 
1-u
\right) 
^{p -1- (p/s)} 
\cdot 
\\
\cdot
\left[ 
\frac{50}{11} 
(p-1)(s-1) + (2ps-2s -p +N) 
(1-u)^{1/s} 
\right]
\\
\qquad \qquad \qquad 
\qquad \qquad \qquad \qquad \qquad \qquad 
\mbox{if }\ 0\leq u\leq 1 
\\
(2s)^{p-1} 
\left[ 
1-(u-1) 
^{1/s} 
\right] 
^{(p/2)-1} 
\cdot
(u-1)^{p-1- (p/s)} 
\cdot 
\\
\cdot 
\left[ 
-2(s-1)(p-1) + (2ps-2s -p+N) (u-1) ^{1/s} 
\right]
\\
\qquad \qquad \qquad 
\qquad \qquad \qquad \qquad \qquad \qquad
\mbox{if }\ 1\leq u\leq 2
\end{array}
\right.
.
\end{equation*}
Note first that
$f\in C^{\infty } ((0,1)\cup (1,2)) \cap C[0,2)$, since
$s > \frac{p}{p-1}$.
Breaking of (radial) symmetry takes place at the level $u=1$ where
$f(1) = 0$.
Let us inspect several cases in more detail.
\\  
{\bf (i)}
Let $p = 2$ (Laplacian case) and $s >  2$. 
\\
Then we have $f \in C^{1-(2/s)} [0,2]$, 
but $f\not\in C^{0,1} [0,2]$. This is not surprising, since Theorem A tells us that the solutions are radially symmetric when $f\in C^{0,1} $. 
\\
{\bf (ii)} If $p\in (1,2)$ and $s> \frac{p}{p-1} $, or if $p>2$ and $s\in (\frac{p}{p-1}, \frac{p}{p-2})$, then  
we have $f \in C^{p-1 - (p/s)}  [0,2)$, 
but $f\not\in C^{0,1} [0,2)$.
\\
{\bf (iii)} If $p>2 $ and $s\geq \frac{p}{p-2}$, then we have $f\in C^{0,1} [0,2]$. 
\\[0.1cm]
Another approach to symmetry has been introduced by one of the authors.
It is based on a 
rearrangement technique called
{\sl continuous Steiner symmetrization\/} 
%(CStS)
(see \cite{Br1}, \cite{Br2}). The method
 allows to obtain 'local' symmetry properties for weak
(distributional) solutions of the following problem,
$$
{\bf (P) }
\left\{ 
\begin{array}{l}
u\in C^1 (\overline{B_R } ),
\\
u\geq 0, \ u \not\equiv 0, \quad 
-{\mathscr L} u \equiv \, -\mbox{div} \left( \displaystyle{\frac{g(|\nabla u|)}{|\nabla u|}}  \nabla u \right) = f(|x|,u)  \quad \mbox{in } \ B_R ,
\\
 u=0 \quad \mbox{on }\ \partial B_R ,
\end{array}
\right.
$$ 
where $g$ and $f$ satisfy certain properties. Roughly speaking, a function is locally symmetric, if it is radially symmetric and radially decreasing in some annuli and flat elsewhere. 
For a precise definition of local symmetry we refer to the next section (Section 2). Note that the function $u$ in Example 1.1 is locally symmetric.
\\ 
\hspace*{1cm} Our paper consists of two parts. After presenting the local symmetry results of \cite{Br2} in Section 2, we will use the Strong Maximum Principle for the degenerate elliptic operator ${\mathscr L}$ to obtain radial symmetry of the solutions, provided that the right-hand side $f$ satisfies appropriate growth conditions near its zero points (see Theorems 3.5 and 3.6 of Section 3). Note that a similar analysis for the $p$--Laplace operator, ($p>1$), that is, for
$$
g(t)= t^{p-1} , \quad (t\geq 0),
$$
has already been carried out by one of the authors in \cite{Br3}.
\\ 
\hspace*{1cm} Still one would like to know whether or not solutions of {\bf (P)} with symmetry breaking (as they appear in Example 1.1) are 'physically relevant'. To this aim we study problem {\bf (P)} and its relation to an associated variational problem in a second part of our paper. For simplicity we restrict ourselves to the case that $f$ is independent of $|x|$. Then solutions of {\bf (P)} are critical points to the variational problem 
\begin{eqnarray*}
J(v) & := & \int_{B_R } \left[ G(|\nabla v|) -F(v) \right] \, dx \longrightarrow \mbox{Inf!} , 
\quad  v\in X,
\\
 & & \mbox{where }\ X:= \{ v\in  
C^1 (\overline{B_R }): \  v=0 \  \mbox{ on }\  \partial B_R \},
\end{eqnarray*} 
$G(t):= \int _0 ^t g(s)\, ds $ and $F(t):= \int_0 ^t f(s)\, ds $, ($t\geq 0$).
Using the local symmetry result of \cite{Br2}, a Pohozaev-type identity, as well as calculating second variations of $J$ in certain directions, we prove that  
local and global minimizers of $J$ in $K$ are radially symmetric under rather mild conditions (see Theorems 4.2, 4.5 and 4.9 of Section 4). 
\section*{2. Local symmetry} 
\setcounter{section}{2} \setcounter{equation}{0} 
In this and the following section we study symmetry properties of solutions of problem {\bf (P)}. We will always assume that the functions $g$ and $f$ satisfy the following properties,
\begin{eqnarray}
\label{gcond1}
 & & g\in C[0, +\infty )\cap C^1 (0,+\infty ),
\\
\label{gcond2}
 & & g(0)=0, \quad g'(t) >0 \quad \mbox{for }\ t>0 ,
\\
\label{fcond1}
 & & f\in C([0,R ]\times [0, +\infty )) , \ f= f(r,t),
\\
\label{fcond2}
 & & \mbox{the mapping }\ r\longmapsto f(r,t) \ \mbox{ is non-increasing }.
\end{eqnarray}
Let us first recall a local symmetry result of \cite{Br2}. 
\\[0.1cm] 
{\bf Definition 2.1.} {\sl Let $u\in C^1 (\overline{B_R }) $ be a nonnegative function with $u\not\equiv 0$. We say that $u$ is locally symmetric if
\begin{eqnarray}
\label{decomp1}
 & & B_R=\bigcup_{k=1}^m A_k \cup \{ x:\nabla u(x)=0\}  ,\quad \mbox{ where} \\
\label{decomp2}
 & & A_k =B_{R_k } (z_k )\setminus \overline{B_{r_k } (z_k )} ,\quad 
(z_k \in B_R ,\ 0\leq r_k 
<R_k ),
\\
\label{decomp3}
 & & u(x) =U_k (|x-z_k |),\quad (x\in A_k ), \ \mbox{ where $U_k \in C^1 [r_k , R_k] $ and } 
\\
\label{decomp4}
 & & U_k ' (r) <0 \quad \mbox{ for } \ r\in (r_k , R_k ), 
\\
\label{decomp5}
 & & u(x)\geq U_k (r_k ) \quad \forall 
x\in B_{r_k } (z_k ),
\\
 & & \qquad \qquad (k=1,\ldots ,m),
\nonumber 
\end{eqnarray} 
the sets $A_k $ are pairwise disjoint 
and $m\in \mathbb{N} \cup \{ +\infty \} $.} 
\\[0.1cm]
{\bf Remark 2.2.} 
The conditions (\ref{decomp1})--(\ref{decomp5}) imply that $u$ is radially symmetric and radially decreasing in annuli $A_k $, $(k= 1, \ldots , m)$, and flat elsewhere in $B_R $.
Note also that, since $u\in C^1 (
\overline{B_R})$, we have that 
\begin{equation}
\label{flatboundary}
U_k ' (r_k )=0, 
\end{equation}
and if $R_k <R$, then also
\begin{equation}
\label{flatboundary2}
U_k ' (R_k ) =0, 
\end{equation}
($  k\in \{ 1,\ldots ,m\} )$.
\\[0.1cm]
{\bf Lemma 2.3.} {\sl (see \cite{Br2}, Theorem 7.2)
\\
Let $u$ be a weak solution of {\bf (P)},
% with $u\not\equiv 0$,  
where the functions $g$ and $f$ satisfy the conditions (\ref{gcond1})--(\ref{fcond2}).   
Then $u$ is locally symmetric.}
\\[0.1cm]
\hspace*{1cm}
In the sequel, we will say that a function $u\in C^1 (\overline{B_R})$ is {\sl radially symmetric and radially non-increasing}, if there is a non-increasing function $U\in C^1 [0,R] $ such that 
\begin{equation}
\label{rad1}
u(x) = U(|x|) \quad \forall x\in B_R .
\end{equation}
{\bf Remark 2.4.} 
Assume that $u$ is as in Lemma 2.3, and that the mapping
$\,  r\longmapsto f(r,t)\, $  
is {\sl strictly } decreasing for $r=|x| $ with $x\in A_k $ and $t\in [U_k (R_k) , U_k (r_k )] $, for some $k$. Then the radial symmetry of $u$ in the annulus $A_k$ implies that $z_k =0$.  Vice versa, if $z_k \not= 0$, then $f$ must be independent of $|x|$ in $A_k $. 
This leads to the following symmetry result.
\\[0.1cm]   
{\bf Corollary 2.5.} {\sl Let $u$, $g$ and $f$ be as in Lemma  2.3, and assume that the mapping 
$\, 
r\longmapsto f(r,t) 
\, $
is strictly decreasing for $r\in [0, R] $ and $t\in [0, \max _{B_R} u ]$. Then $u$ is radially symmetric and radially non-increasing.}
\section*{3. Radial symmetry}
\setcounter{section}{3} 
\setcounter{equation}{0} 
\hspace*{0.7cm}In this section we will use the Strong Maximum Principle for
the operator ${\mathscr L}$ to obtain radial symmetry of
the solutions of problem {\bf (P)} under some additional conditions on $f$.
Our results and proofs are modelled after \cite{Br3} where
a similar  analysis was carried out for the $p$-Laplacian.  
\\
\hspace*{1cm} First we recall a general version of the {\sl Strong Maximum Principle}, see  \cite{PuSeMP}.
\\
Let $g\in C[0,+\infty )$, $g $ strictly increasing and $g(0)=0$.
Further, assume $\varphi \in C[0,+\infty )$, $\varphi (t) \geq 0 $ for all $t$, and  $\varphi (0)=0$.
Finally, let $u\in C^1 (\Omega )$ be a distributional solution to 
\begin{equation}
\label{mainineq}
{\mathscr L} u \equiv \, \mbox{div} \left( \displaystyle{\frac{g(|\nabla u|)}{|\nabla u|}} \nabla u \right) \leq \varphi (u), \quad u\geq 0 \quad \mbox{in }\ \Omega  ,
\end{equation}
where $\Omega $ is a domain in $\mathbb{R}^N $.
Here and in the following we make the convention that
\begin{equation*}
   \frac{g(|y|)}{|y|} y= 0 \, ,
  \quad\mbox{ if $y = 0$ in $\mathbb{R}^N$. }
\end{equation*}
By the Strong Maximum Principle (SMP) for (\ref{mainineq}) to hold, we mean the statement that 
{\sl if $u$ is a 
solution of (\ref{mainineq}) with $u(x_0 )=0$ for some $x_0 \in \Omega $, then $u\equiv 0$ in $\Omega $.}
\\[0.1cm]
It will be convenient to work with the following definition.
\\[0.1cm]
{\bf Definition 3.1.}
{\sl A function $\varphi : [0,+\infty ) \to [0, +\infty )$ belongs to
the class ${\mathscr A}_g$
(i.e., $\varphi \in {\mathscr A}_g$),
if $\varphi\in C[0,\infty )$, $\varphi (0) = 0$ and either
$\varphi (t) = 0$ for $t\in [0,d]$, ($d > 0$), or else
$\varphi (t) > 0$ for $t\in (0,\delta)$, ($\delta > 0$), and
\begin{equation}
\label{degenerateint}
  \int_0^{\delta} \frac{dt}{ H^{-1}\left( \Phi(t)\right) } = \infty \,,
\end{equation}
where the functions $H$ and $\Phi$ are given by
$$
H(t):= tg(t) - \int_0 ^t g(s)\, ds , \quad \Phi (t):= \int_0 ^t \varphi (s)\, ds , \quad (t\geq 0).
$$
}
{\bf Lemma 3.2. } {\sl ( Strong Maximum Principle (SMP), see 
\cite{PuSeMP}, Theorem 1.1.1.)
\\
 In order for the strong maximum principle to hold for 
(\ref{mainineq}) 
it is necessary 
and sufficient that $\varphi \in {\mathscr A}_g $.
}
\\[0.1cm]
{\bf Remark 3.3.}
In case of the $p$-Laplace operator we have
$g(t) = t^{p-1}$, ($p>1$), and condition \eqref{degenerateint} reads
\begin{equation}
\label{degenerateintp}
\int_0 ^{\delta } \frac{dt}{ [\Phi (t)] ^{1/p} } = \infty .
\end{equation}
In particular, (\ref{degenerateintp}) holds if
\begin{equation}
\label{phiplaplace}
\varphi (t) \leq  Ct^{p-1} , \quad (t\in (0, \delta )),
\end{equation}
for some $C>0$.
\\[0.1cm]
It is well-known that the SMP implies the following 
\\[0.1cm]
{\bf Lemma 3.4.}
{\sl (Boundary Point Lemma, see \cite{PuSeMP}, Theorem 5.5.1.)
\\ 
Let $\Omega $, $g$, $\varphi$, and $u$ be as in Lemma 3.2,
with $u>0$ in $\Omega$,  and let $y$ be a point on $\partial\Omega$
satisfying the interior sphere condition with a ball $B\subset \Omega$.
Let $\nu = \nu(y)\in \mathbb{R}^N$ be
the corresponding exterior normal at
$y\in \partial\Omega\cap \partial B$.
Then the exterior normal derivative of $u$ at $y$ satisfies
$$
\frac{\partial u}{\partial \nu }(y) < 0 \quad\mbox{ at } y \,.
$$
}
\hspace*{1cm}Before stating the main results of this section, we make a simple observation.
\\[0.1cm]
{\bf Lemma 3.5.} $\ $ {\sl Let $u$ be as in Lemma 2.3, and let 
$A = B_{\rho ' } (z) \setminus \overline{B_{\rho}(z)}$
be one of the annuli in the decomposition \eqref{decomp1},  with $z\in B_R $, $0\leq \rho <\rho' <R$  and
$x_0\in \partial A\cap B_R$.
Then, if 
$x_0\in \partial B_{\rho}(z) $ and $\rho > 0$, or if 
$x_0\in \partial B_{\rho'}(z)$ and $N\geq 2$, we have that}
\begin{equation}
\label{f=0}
f(|x_0 |,u(x_0 ))=0.
\end{equation}
{\sl Proof :\/} $\ $
First assume that $x_0\in \partial B_{\rho} (z)$ and $\rho > 0$.
By (\ref{decomp5}) and (\ref{flatboundary}) we have that 
$u(x)\geq u(x_0)$ in $B_{\rho} (z)$ and 
\begin{equation}
\label{nab0} 
\nabla u(x_0) = 0.
\end{equation}
Assume that $f(|x_0|,u(x_0)) > 0$.
Then, by the continuity of the functions $f$ and $u$,
there is a number $\delta > 0$ such that
$f(|x|,u(x)) > 0$ holds for all
$x\in B_R\cap B_{\delta}(x_0)$.
Since $u$ is a solution of {\bf (P)}, we may apply Lemma 3.2 to conclude that $u(x) >u(x_0 )$ for all $x\in B_{\rho } (z )\cap B_{\delta } (x_0 ) $. Then,  
Lemma 3.4 tells us that
\begin{equation}
\label{hopf1}
\frac{\partial u}{\partial \nu }(x_0) < 0 \ \mbox{ at }\, x_0 , \quad (\nu :\ \mbox{ exterior normal to $B_{\rho } (z)$}).
\end{equation}
But this  contradicts (\ref{nab0}).
\\
We have verified that $f(|x_0|,u(x_0))\leq 0$.
On the other hand, since $u(x)\leq u(x_0)$ in $A$,
Lemma 3.4 also gives
$f(|x_0|,u(x_0))\geq 0$, by a similar reasoning.  This proves (\ref{f=0}) in the considered case.  
\\
Next assume that $x_0\in \partial B_{\rho'}(z)$ and $N\geq 2$.
Since
$u(x)\geq u(x_0)$ in $A$ and $\nabla u(x_0) = 0$,
Lemmata 3.3 and 3.4 again yield $f(|x_0|,u(x_0))\leq 0$.
It remains to show that
\begin{equation}
\label{fpositive}
  f(|x_0|,u(x_0))\geq 0.
\end{equation}
Assume  that 
\begin{eqnarray}
\label{limC}
 & & \mbox{there is a subsequence $\{ A_{k'} \} $ of annuli in (\ref{decomp1}), such that}
\\
\nonumber  
 & &  z_{k'} \to x_0 \ \mbox{ and } \ R_{k'} \to 0 \ \mbox{ as }\  k'\to \infty .
\end{eqnarray}
Since 
$u(x) > u\Big|_{\partial B_{R_{k'}}(z_{k'})} $ in $B_{R_{k'}}(z_{k'})$,
the Strong Maximum Principle tells us that there are points
$y_{k'} \in B_{R_{k'}} (z_{k'}) $ such that 
$0\leq -{\mathscr L}  u(y_{k'}) =f(|y_{k'} | ,u(y_{k'}))$.
Since also $ y_{k'} \to x_0 $, we obtain (\ref{fpositive}) in this case.\\
Now suppose that (\ref{limC}) does not hold. Then one of the following situations {\bf (i)} or {\bf (ii)} occurs: 
\begin{eqnarray*}
 {\bf (i)} \qquad & & \mbox{there is an annulus in the decomposition (\ref{decomp1}), denoted by  }
\\
 & & \mbox{$A_1 = B_{R_1} (z_1) \setminus \overline{B_{r_1 }
(z_1 ) } $, with $x_0 \in \partial B_{R_1 } (z_1 )$  and}
\\
 & & \mbox{$ \nabla u(x)=0$ in $ \ B_{\varepsilon } (x_0 ) \setminus \overline{(A\cup A_1 )}$ 
for some $ \varepsilon >0;\ $  or}
\\ 
 {\bf (ii)} \qquad & &  \nabla u(x)=0 \ \mbox{ in } \ 
B_{\varepsilon } (x_0 ) \setminus \overline{A}
\ \mbox{\ for some } \ \varepsilon >0.
\end{eqnarray*}
 Clearly in both cases {\bf (i)} and {\bf (ii)} we have 
$u(x)=u(x_0 ) $  on some open subset of $ B_{\varepsilon} (x_0 ) $,
which means that $f(|x_0 |,u(x_0 ))=0$. 
$\hfill \Box $
\\[0.1cm]
{\bf Remark 3.6.} The last step of the above proof does not work in the case $N=1$, because  we cannot deduce that $u$ is constant on an open subset of $B_{\varepsilon } (x_0 )$. 
\\[0.1cm]
\hspace*{1cm}Lemma 3.5 allows to exclude symmetry breaking for solutions of problem {\bf (P)} if $f$ satisfies appropriate growth conditions near its zeros.  
\\[0.1cm] 
{\bf Lemma 3.7.}   
{\sl Let $u$ be a solution of problem {\bf (P)} and suppose that $f$ satisfies the following condition:
\begin{eqnarray*}
%\nonumber
{\bf (a)}
 & & \mbox{ if $f(\rho ,\tau )=0$ for some $\tau > 0$ and $\rho\in [0,R)$,}
\\
 & & \mbox{  
then there is a number $\delta \in (0, \tau )$ and a function $\varphi\in {\mathscr A}_g$  such that}
\\
 & & \mbox{ 
$f(r,t)\leq \varphi (\tau -t)$  for   $t \in [\tau - \delta , \tau ] $  and  
$ r \in [\rho -\delta , \rho + \delta ]$.}
\end{eqnarray*}
Then there are mutually disjoint balls $B_{R_k } (z_k )  $, ($z_k \in B_R $, $R_k >0$), such that 
\begin{eqnarray}
\label{decomp6}
 & &  B_R =\bigcup_{k=1}^m B_{R_k } (z_k ) \cup \{ x:\, u(x)=0 \} ,
\\
\label{decomp7}
 & & u(x)= U_k (|x-z_k |), \quad (x\in B_{R_k } (z_k )), \mbox{ where $U_k \in C^1 [0,R_k ] $, }
\\
\label{decomp8}
 & & 
\mbox{
$U_k (R_k )= 0$, $ U_k ' (r) <0$ for $\ r\in (0, R_k )$, 
}
\\
 \nonumber & & \quad \qquad (k=1, \ldots , m ),
\end{eqnarray}
and $m\in \mathbb{N} \cup \{ +\infty \} $. }
\\[0.1cm]
{\sl Proof:\/} $\ $
By Lemma 2.3, $u$ is locally symmetric. 
First we claim that the annuli $A_k $ in (\ref{decomp2}) 
are in fact {\sl punctured balls\/},
i.e., we have 
\begin{equation}
\label{rk=0}
r_k =0 ,\quad (k=1,\ldots ,m).
\end{equation}
Let $A=B_{\rho '} (z)\setminus \overline{B_{\rho } (z)} $, ($0\leq \rho <\rho' \leq R $), 
be one of the annuli in (\ref{decomp2}) and suppose that $\rho >0$.
By Lemma 3.5., we have
$f(|x|,u_0 )=0 \ $ $\ \forall x\in \partial B_{\rho} (z)$, where 
$ u_0 =u \Big|_{\partial B_{\rho } (z)} $.
We set $w:= u_0 -u$. Then
$w\in C^1 (\overline{A}) $, $w> 0$ in $A$ and $w= |\nabla w|=0$ on $\partial B_{\rho } (z)$. Furthermore, assumption {\bf (a)} yields 
$$
{\mathscr L}w = f(|x|,u_0 -w) \leq \varphi(w)
\quad \mbox{ in }\, A \,,
$$
for some function $\varphi\in {\mathscr A}_g$.
But this is impossible, by Lemma 3.4.
Hence, we must have $\rho  = 0$.
This proves \eqref{rk=0}.
\\
Next we claim that  
\begin{equation}
\label{u=0}
  u = 0 \quad\mbox{ on }\; \partial B_{R_k}(z_k) \,;
        \quad k=1,2,\ldots,m \,.
\end{equation}
We fix some $l\in \{ 1,2,\ldots,m\}$, and set
\begin{equation*}
  u_l (x) := \left\{ 
\begin{array}{ll}
  u\Big|_{\partial B_{R_k}(z_k)}
& \mbox{ if }\; x\in B_{R_k}(z_k) \,;
  \quad k=1,2,\ldots,m \,,\; k\not= l \,,
\\
  u(x) & \mbox{ otherwise. }
\end{array}
\right.
\end{equation*}
Since $\nabla u=0$ on $\partial B_{R_k } (z_k )$ whenever $R_k <R $, ($k= 1, \ldots ,m$),  we have that $u_l \in C^1 (\overline{B_R })$. 
Hence
$\nabla u_l =0$ in $B_R \setminus B_{R_l } (z_l )$ which means that 
$u_l (x) =u\Big| _{\partial B_{R_l } (z_l )}  $ for all $x\in 
B_R \setminus 
B_{R_l } (z_l )$. 
Since $u_l \leq u$, this implies that $u=0$ on $\partial B_{R_l } (z_l )$, 
and (\ref{u=0}) follows. Now the assertions of the Lemma follow from (\ref{rk=0}) and (\ref{u=0}).
$\hfill \Box $
\\[0.1cm]
{\bf Lemma 3.8.} 
{\sl  Let $u$ be a solution of problem {\bf (P)} where $N=2$, and suppose that $f$ satisfies the following condition,
\begin{eqnarray*}
{\bf (b)} & & 
\mbox{ if $f(\rho ,0)=0$ for some $\rho \in [0,R )$, 
then there is a number $\delta >0$ and a function}
\\
 & & \mbox{ $\varphi \in {\mathscr A} _g $, 
such that $ -f(r,t)  \leq \varphi (t) $ for $\ r\in [0, R)$ and $t \in [0, \delta ] $.} 
\end{eqnarray*}
Then $u>0 $ in $B_R $.}
\\[0.1cm] 
{\sl Proof: } Suppose that $u $ is not positive in $B_R $. Setting $\Omega :=\{ x:\, u(x)>0\}$, we choose a point $x_0 \in \partial \Omega \cap B_R  $ which satisfies an interior sphere condition. Then $u(x_0 )=0$ and $\nabla u(x_0 )=0$. Moreover, Lemma 2.3 tells us that $x_0 $ is also a boundary point of one of the annuli in the decomposition (\ref{decomp1}). By Lemma 3.5 this implies $f(|x_0 |,0)=0$. 
Hence condition {\bf (b)} is satisfied (with $|x_0 |=\rho $). Then Lemma 3.4 yields $|\nabla u(x_0 )|>0$, a contradiction. 
$\mbox{ } $ 
$\hfill \Box $
\\[0.1cm]
{\bf Theorem 3.9.} 
 {\sl Let $u$ be a solution of problem {\bf (P)}, and suppose that $f$ satisfies the following condition,
\begin{eqnarray*}
{\bf (c)}
 & & \mbox{ if $f(\rho ,\tau )=0$ for some $\tau \geq 0  $ and $\rho\in [0,R)$,}
\\
 & & \mbox{  
then there is a number $\delta >0$ and a function $\varphi\in {\mathscr A}_g$  such that}
\\
 & & \mbox{ 
$-f(r,t)\leq \varphi (t-\tau )$  for   $\, t \in [\tau , \tau +\delta ] $  and  
$\,  r \in [\rho -\delta , \rho + \delta ]$.}
\end{eqnarray*}
Then $u$ is radially symmetric and radially non-increasing. Moreover, there is a number $R' \in [0,R)$ such that
\begin{eqnarray}
\label{rad2}
 & & U'(r) <0 \ \mbox{ for $\, r\in (R', R)$, and }
\\
\label{rad3}
 & & U(r) = \mbox{ const. } \mbox{ for $\,  r\in [0,R']$.}
\end{eqnarray}
}
{\sl Proof: } 
Let $A= B_{\rho ' } (z) \setminus \overline{B_{\rho } (z)}$ be one of the annuli in the decomposition (\ref{decomp1}) with $\rho >0$, and let $u_0 := u|_{\partial B_{\rho } (z)} $. We claim that
\begin{equation}
\label{uconst}
u(x)= u_0 \ \mbox{ in $B_{\rho } (z)$.}
\end{equation}
Assume that this is not true, and set $\Omega := \{ u>u_0 \} \cap B_{\rho } (z)$. Let $x_0  $ be a boundary point of $\Omega $ satisyfying an interior sphere condition. Setting $w:= u-u_0 $, we have $w(x_0 )= |\nabla w(x_0 )|=0$. By assumption {\bf (c)} we have that 
$$
\mathscr{L} w= -f(|x|, w+u_0 ) \leq \varphi (w) 
$$
in a neighborhood of $x_0 $, while $w>0 $ in $\Omega $. By Lemma 3.4 this implies that $|\nabla w (x_0 )|>0$, a contradiction. 
\\
Further, by Lemma 3.8 we have that $u>0$ in $B_R $. Together with (\ref{uconst}) this shows that $m=1$ in (\ref{decomp1}) and $A_1 = B_R \setminus \overline{B_{R'} }$, for some $R' \in [0, R)$. The Theorem is proved.     
\hfill \Box
\\[0.1cm]
{\bf Theorem 3.10.}   
{\sl Let $u$ be a solution of problem {\bf (P)} and suppose that conditions {\bf (a)} and {\bf (b)} are satisfied. If $N=1$, we additionally assume that 
\begin{equation}
\label{f>0}
f(r,0) \geq 0 \quad \mbox{for $\ r\in [0,R]$.}
\end{equation}
Then $u$ is radially symmetric and radially non-increasing. Moreover, there holds
\begin{equation}
\label{RadU}
U'(r) <0 \quad \mbox{for $\ r\in (0,R)$.}
\end{equation}
} 
{\sl Proof:\/} $\ $
By Lemma 3.7, $u$ satisfies (\ref{decomp6})-(\ref{decomp8}). We split into two cases.
\\
{\bf 1. }  Let  $N\geq 2$. Then we have that $u>0$ by Lemma 3.8. This implies that $m=1$,  $z_1 =0$ and $R_1 =R$, so that (\ref{RadU}) follows.
\\
{\bf 2. } Let $N=1$. Assume that 
for some $k\in \{ 1,2,\ldots,m\}$,
\begin{equation}
\label{Rk<R}
  B_{R_k}(z_k)\not= B_R \,.
\end{equation}
Then we have 
$u=  |\nabla u| = 0$  on $ \partial B_{R_k}(z_k) $. Then Lemma 
3.2 yields 
$f(|z_k \pm  R_k|,0)\leq 0$, while  assumption 
(\ref{f>0}) rules out the cases $f(|z_k +  R_k|,0) < 0$ and $f(|z_k - R_k|,0)< 0$.
Hence, we must have $f(|z_k \pm R_k |,0) = 0$. Then, 
proceeding as in the proof of Lemma 3.8, we deduce that $u>0$ in $B_R $, which again implies $m=1$, $z_1 =0$ and $R_1 =R$.  
$\hfill\Box$
\\[0.1cm]
{\bf Remark 3.11.} Let us specify some typical situations where the above conditions {\bf (a)}, {\bf (b)} and {\bf (c)} are satisfied.  
\\
{\bf 1. } If $f(r,t)\geq 0$, ($r\in [0,R]$, $t\geq 0$), then conditions {\bf (b)} and {\bf (c)} hold.
\\
{\bf 2. } If $f(r,t)>0$, or if $t\longmapsto f(r,t)$ is non-decreasing, ($r\in [0,R]$, $t\geq 0$), then the conditions {\bf (a)}, {\bf (b)} and {\bf (c)} hold. 
\\
{\bf 3. } If $g(t)= t^{p-1}$, ($p$-Laplace operator), $p\in (1,2]$, and if 
$$
f= f_1 + f_2,
$$
where $f_i \in C([0,R]\times [0, +\infty ))$, $r\longmapsto f(r,t)$, nonincreasing, ($i=1,2$), $f_1 (r, \cdot )$ is H\"older continuous with exponent $p-1$, uniformly for all $r\in [0,R]$, and $t\longmapsto f_2 (r,t) $ is non-decreasing, ($r\in [0,R]$, $t\in [0,+\infty )$),  then the conditions {\bf (a)}, {\bf (b)} and {\bf (c)} hold. 
\section*{3. Instability of non-symmetric solutions}
\renewcommand{\theequation}{\thesection.\arabic{equation}}
\setcounter{section}{4} \setcounter{equation}{0} 
In this section we study radially symmetric solutions in an annulus.
We focus on their stability w.r.t. the associated variational problem.
\\[0.1cm]
We assume 
\begin{eqnarray}
\label{assg}
 & & g\in C[0,+\infty )\cap C^1(0,+\infty) \,,\ g(0) = 0 \,,\
     \mbox{ and $g$ is strictly increasing; }
\\
\label{assf}
 & & f\in C[0,+\infty) \,.
\end{eqnarray}
We set
\begin{equation*}
\textstyle
  G(t) := \int_0^t g(s)\, ds \,,\quad
  F(t) := \int_0^t f(s)\, ds \,;\quad t\in [0,+\infty) \,.
\end{equation*}
Let $\Omega$ be a bounded domain in $\mathbb{R}^N$ with $C^1$-boundary.
\\
We consider 
weak solutions to the following boundary value problem,
\begin{equation} 
\label{bvp}
\left\{ 
\begin{array}{l} 
u\in C^1 (\overline{\Omega }) , \\
 - \mbox{div} \left( \displaystyle{\frac{g(|\nabla u| )}{|\nabla u|}} \nabla u \right) = f(u) \ 
\mbox{ in } \Omega ,
\\
 u=0 \ \mbox{ on } \partial \Omega .
\end{array}
\right.
\end{equation}
The solutions $u$ are critical points to the variational functional
$$
  J(v) := \int_{\Omega } \left( G(|\nabla v|) -F(v) \right)
          \, dx \,,\quad v\in X \,, 
$$
where
$$
  X:= \{ v\in C^1(\overline{\Omega}) :\,
         v = 0 \;\mbox{ on } \partial\Omega \} \,.
$$
We say that $u$ is
a {\sl global minimizer for $J$ in $X$} if $u\in X$ and
$$
  J(u)\leq J(v) \quad\mbox{ holds for all }\; v\in X \,.
$$
Furthermore, we say that $u$ is
a {\sl local minimizer for $J$ on $X$} if $u\in X$ and
there exists an $\varepsilon > 0$ such that
$$
  J(u)\leq J(v) \quad\mbox{ whenever $v\in X$ satisfies }\;
  \Vert v-u\Vert_{C^1(\Omega)} < \varepsilon \,.
$$
Our first result treats the Neumann boundary conditions:
\\[0.1cm]
{\bf Lemma 4.1.}
{\sl Let $u$ be a global minimizer to $J$ in $X$, and assume that 
\begin{equation}
\label{unu=0}
\frac{\partial u}{\partial \nu } =0 \ \mbox{ on } \partial \Omega .
\end{equation}
Then we have $u\equiv 0 $ in $\Omega $.
}
\\[0.1cm]
{\sl Proof :}
Since $u$ is a critical point for $J$ in $X$, it is a solution of
the boundary value problem \eqref{bvp}.
Hence, for every vector field
\begin{math}
  h\in C^1\left( \overline{\Omega}, \mathbb{R}^N\right) ,
\end{math}
the following integral identity of Pohozhaev\--type holds
(see \cite{PucciSerrin}, formula~(4)),
\begin{eqnarray}
\label{PucSer1}
 & & 
\int_{\partial \Omega } \Big\{ G(|\nabla u|) -g(|\nabla u|) |\nabla u| \Big\} \, (h\cdot \nu )\, {\mathscr H}_{N-1} (dx) 
\\
 \nonumber
 & = & \int_{\Omega } \left\{ \left[ G(|\nabla u|) -F(u)\right] \, \mbox{div}\, h - \sum_{i,j=1}^N (h_j)_{x_i }  u_{x_j} u_{x_i } \frac{g(|\nabla u|)}{ |\nabla u|} \right\} \, dx .
\end{eqnarray}
Choosing $h(x) =x$ this becomes
\begin{eqnarray}
\label{PucSer2}   
 & & 
\int_{\partial \Omega } \Big\{ G(|\nabla u| )- g(|\nabla u|) |\nabla u| \Big\} \, (x\cdot \nu )\, {\mathscr H}_{N-1} (dx) 
\\
 \nonumber
 & = & \int_{\Omega } \Big\{ N[G(|\nabla u|) -F(u)]  - g(|\nabla u|) |\nabla u| \Big\} \, dx .
\end{eqnarray}
The left\--hand side of \eqref{PucSer2} is zero by \eqref{unu=0}.
Hence, we further deduce
\begin{equation}
\label{PucSer3}   
  N J(u) =
  N \int_{\Omega } \Big\{ G(|\nabla u|) - F(u)\Big\} \,\mathrm{d}x
  = \int_{\Omega } g(|\nabla u|) |\nabla u| \,\mathrm{d}x \,.
\end{equation}       
The right\--hand side of \eqref{PucSer3} is positive unless $u\equiv 0$,
in which case we have $J(u) = 0$.
The lemma is proved.
$\hfill\Box$
\\[0.1cm] 
Now we show that if $\Omega $ is a ball, then nonnegative global minimizers of $J$ in $X$ are radially symmetric.
\\[0.1cm]
{\bf Theorem 4.2.}
{\sl Let $u$ be a solution of {\bf (P)} which is also a global minimizer for $J$ in $X$, with $\Omega =B_R$.
Then $u$ is positive in $B_R$,
radially symmetric and radially non-increasing.
Moreover, conditions (\ref{rad2}) and (\ref{rad3}) hold.
}
\\[0.1cm]
{\sl Proof :} $u$ is locally symmetric by Lemma 2.3.
Suppose that $u$ is not radially symmetric.
Then there exists a ball $B := B_{R'}(x_0)\subset B_R$, 
($x_0\in B_R$, $R'\in (0,R)$), such that
$\frac{\partial u}{\partial\nu} = 0$,
$u(x) =: u_0\in [0,+\infty)$ on $\partial B$, and
$u(x)\geq u_0$, $u(x)\not\equiv u_0$ in~$B$.
We define
\begin{eqnarray*}
  v(x) & := &
\left\{
\begin{array}{ll} 
  u(x) &\mbox{ if }\; x\in B_R\setminus B \,,
\\
  u_0  &\mbox{ if }\; x\in\overline{B} \,;
\end{array}
\right.
\\
  \overline{u}(x) & := &
  u(x) - u_0 \quad\mbox{ for $x\in \overline{B}$. }
\end{eqnarray*}
Then $v\in X$ and 
\begin{eqnarray}
\label{ineq1}
  0 &\geq & J(u) - J(v)
  = \int_{B}\Big( G(|\nabla u|) - F(u) + F(u_0)\Big) \,\mathrm{d}x
\\
\nonumber
 & = & \int_{B}
  \Big( G(|\nabla\overline{u}|) - F(\overline{u} + u_0) + F(u_0)
  \Big) \,\mathrm{d}x \,.
\end{eqnarray}    
Furthermore, we have
$\overline{u}\in C^1(\overline{B})$,
\begin{math}
  \overline{u} = \frac{\partial\overline{u}}{\partial\nu} = 0
\end{math}
on $\partial B$, together with $\overline{u}\geq 0$ and
${}- {\mathscr L}\overline{u} = f(\overline{u} + u_0)$ in 
$B$.
Hence, $\overline{u}$ is a critical point for the variational functional
\begin{equation*}
  \overline{J}(w) := \int_B
  \Big( G(|\nabla w|) - F(w + u_0) + F(u_0) \Big) \,\mathrm{d}x
  \quad\mbox{ for }\; w\in \overline{X} \,,
\end{equation*}
where
\begin{math}
  \overline{X} := \{ w\in C^1(\overline{B}):\, w = 0
                     \,\mbox{ on }\, \partial B\} .
\end{math}
Using \eqref{ineq1}, this yields
$\overline{J}(\overline{u})\leq 0$.
On the other hand, repeating the calculation of the previous proof
-- with $B$ and $(F(w + u_0) - F(u_0))$ in place of $B_R$ and $F(w)$,
respectively -- we obtain
$\overline{J}(\overline{u})\geq 0$.
In turn, in analogy with \eqref{PucSer3}, this yields
\begin{math}
    N\overline{J}(\overline{u})
  = \int_B g(|\nabla\overline{u}|)\, |\nabla\overline{u}|\,\mathrm{d}x
  = 0
\end{math}
and, thus, $\overline{u}\equiv 0$.
Hence, $u\equiv u_0$ in $B$, a contradiction.
$\hfill\Box$
\\[0.1cm]
\hspace*{1cm} Next, our aim is to show that also nonnegative {\sl local} minimizers to $J$ in a ball are radially symmetric. 
As a first step, we examine local minimizers in starshaped domains.
\\
We call a domain $\Omega $ {\sl starshaped  w.r.t. the origin } if for any $x\in \Omega $ we have $\{tx:\, 0\leq t\leq 1 \} \subset \Omega $. 
\\[0.1cm]
{\bf Lemma 4.3.} {\sl 
Let $\Omega $ be starshaped w.r.t. the origin.
Furthermore, let $u$ be a solution to \eqref{bvp} satisfying \eqref{unu=0},
such that $u$ is a local minimizer for $J$ in~$X$.
Finally, assume that 
\begin{equation}
\label{mon1}
  \mbox{ the mapping }\; t\longmapsto g(t)\, t^{1-N}\;
  \mbox{ is strictly decreasing. }
\end{equation}
Then $u\equiv 0$ in $\Omega $.}
\\[0.1cm]
{\bf Remark 4.4. }
Condition \eqref{mon1} is satisfied, for instance, if 
$g(t)= t^{p-1} $ (the $p$-Laplace operator) with $1 < p < N$.
\\[0.1cm]
{\sl Proof of Lemma 4.3: } 
First observe that condition \eqref{mon1} is equivalent to
$$
  g(t(1+s)) - (1+s)^{N-1} g(t) < 0
  \quad\mbox{ for all $s>0$ and $t>0$. }
$$
In turn, recalling
$G'= g$ and $F'= f$ on $\mathbb{R}_+ = [0,+\infty)$,
this implies that the mapping
\begin{equation}
\label{mon2}
  s\longmapsto G(t(1+s)) - \frac{1}{N} (1+s)^N g(t)\, t
  \quad\mbox{ is strictly decreasing for all $t>0$. }
\end{equation}
Now we set for any $\varepsilon > 0$, 
$$
  u_{\varepsilon} (x) :=
\left\{ 
\begin{array}{ll}
  u((1+\varepsilon)x) &\mbox{ if }\; (1+\varepsilon)x\in \Omega  \,,
\\
  0 &\mbox{ if }\; (1+\varepsilon)x\not\in \Omega \,.
\end{array}
\right.
$$
Then $u_{\varepsilon} \in X$.
In view of \eqref{PucSer3} and \eqref{mon2} we have
\begin{eqnarray}
\nonumber
J(u_{\varepsilon }) & = &
\int _{\Omega } \Big( G(|\nabla u| (1+\varepsilon )) - F(u )\Big) (1+\varepsilon )^{-N} \, dx
\\
\nonumber 
 & = & 
\int _{\Omega } \Big(  G(|\nabla u| (1+\varepsilon )) - G(|\nabla u|)  +\frac{1}{N} g(|\nabla u|) |\nabla u| \Big) (1+\varepsilon )^{-N} \, dx
\\
\label{J(ueps)}
 & \leq & \frac{1}{N} \int_{\Omega } g(|\nabla u|) |\nabla u|\, dx =J(u),
\end{eqnarray}
with equality only if $\nabla u = 0$ in $\Omega$.
From this the assertion follows.
$\hfill\Box$
\\[0.1cm]
{\bf Theorem 4.5. } {\sl Let $u$ be a solution of problem {\bf (P)} which is also a local minimizer for $J$ in $X$ with $\Omega = B_R $. Furthermore, assume that (\ref{mon1}) is satisfied. Then $u$ is radially symmetric and radially non-increasing. Moreover, (\ref{rad2}) and (\ref{rad3}) hold.
}
\\[0.1cm]
{\sl Proof: } Suppose there exists a ball $B_{\rho } (z)  \subset \subset B_R $ and a number $u_0 \geq 0 $ such that $u= u_0 $ and $\partial u/\partial \nu =0 $ on $\partial B_{\rho } (z) $ and $u\geq u_0 $ in $B_{\rho} (z)$. For every $\varepsilon >0$ we define 
$$
u_{\varepsilon } (x) := \left\{
\begin{array}{ll}
u(x) & \mbox{ if $x\in B_R \setminus \overline{B_{\rho } (z)}$}
\\
u_0 & \mbox{if $\rho  /(1+\varepsilon)  \leq |x-z| <\rho $}
\\
u (
(1+\varepsilon ) x-\varepsilon z) & \mbox{ if $x\in B_{\rho /(1+\varepsilon)} (z) $ }
\end{array}
\right.
.
$$
Then $u_{\varepsilon } \in X $ and proceeding analogously as in the last proof we obtain 
\begin{eqnarray*}
J(u_{\varepsilon} )- J(u) & = & \int_{B_{\rho } (z)} \Big( G(|\nabla u|(1+ \varepsilon ) -F(u)\Big)  (1+\varepsilon )^{-N} \, dx
\\
 & & - \int_{B_{\rho} (z)} \Big( G(|\nabla u|) - F(u)\Big) \, dx
\\
 & \leq & 0,
\end{eqnarray*}
with equality only if $\nabla u=0$ in $B_{\rho } (z)$.
Since $u$ is a local minimizer, this implies that $u=\mbox{const.}$ in $B_{\rho } (z)$, and the assertion follows.     
$\hfill \Box $
\\[0.1cm] 
\hspace*{1cm}Next we examine {\sl radially symmetric} solutions of  the problem
\begin{eqnarray}
\label{annulusbvp} 
\left\{ 
\begin{array}{l}
u\in C^1 (\overline{A}) ,
\\
- \mbox{div} \left( \displaystyle{\frac{g(|\nabla u| )}{|\nabla u|}} \nabla u \right) = f(u) \ 
\mbox{ in } A ,
\end{array}
\right.
\end{eqnarray}
where $A$ is the annulus
$$
  A := B_{R_2}\setminus \overline{B_{R_1}} , \quad (0 < R_1 < R_2).
$$
{\bf Lemma 4.6.} {\sl
Let $u$ be a distributional solution of (\ref{annulusbvp}) which is radially symmetric. 
Moreover, assume:
\\
{\bf (i)}
There exists a function $U\in C^1(R_1, R_2)$ such that
$u(x) = U(|x|)$ in $C$, $U^{\prime}(r)<0$ for $r\in (R_1, R_2)$,
and $U^{\prime}(R_1) = U^{\prime}(R_2) = 0$;
\\
{\bf (ii)}
$f\in C^1(u_2, u_1)\cap C([u_2, u_1])$, where
$U(R_1) = u_1$ and $U(R_2) = u_2$;
\\
{\bf (iii)}
$f(u_1) = f(u_2) = 0$;
\\
{\bf (iv)}
there exist positive constants $\Gamma$ and $t_0 $ such that
\begin{equation}
\label{technass}
t^2 g'(t) \leq \Gamma g(t) \quad \mbox{ for $0<t\leq t_0 $.}
\end{equation}
Then there exists a function $\Phi\in C_0^1(R_1, R_2 )$
such that}
\begin{equation}
\label{2ndvar<0}
\int_{R_1 } ^{R_2 }  \Big\{  
g^{\prime } (-U^{\prime } ) (\Phi ^{\prime } ) ^2 - f^{\prime } (U) \Phi ^2 \Big\}  
r^{N-1 } \, dr < 0.
\end{equation}     
{\bf Remark 4.7.} 
\\
{\bf (a)}
Notice carefully that we do {\sl not\/} assume that
$f$ is $C^1 $ at $u_1 $ and $u_2 $.
\\
{\bf (b)}
The assumption \eqref{technass}) is technical, but it is essential
in the proof that we present below.
Note also, that it is satisfied for many relevant differential operators.
Here are some examples:
\\ 
{\bf 1. }\ $g(t) = t^{p-1} $, where $p>1$, ($p$-Laplace operator);
\\
{\bf 2. }\ $g(t) = \sum_{k=1}^m c_k\, t^{p_k - 1}$, where 
$c_k >0, \, p_k >1$, ($k=1, \ldots m$);
\\
{\bf 3. }\ $g(t) = \frac{t}{\sqrt{1+t^2}}$, (minimal surface operator);
\\
{\bf 4. }\ $g(t) = e^{- \Gamma t^{-\alpha} }$
          with $\alpha\in (0,1]$ and $\Gamma > 0$.
\\
It is not satisfied, for instance, for
$g(t) = e^{- \Gamma t^{-\alpha}}$ when $\alpha > 1$ and $C>0$.
\\[0.1cm]
{\sl Proof of Lemma 4.6. :}   The idea is to use an appropriate cut--off function of $U^{\prime } $ for $\Phi $.  
\\
We put
$\Phi :=  U^{\prime} H$, where 
$H\in C_0 ^1 (R_1 , R_2 ) $. 
Then we evaluate
\begin{eqnarray}
\label{QH} & &
\\
\nonumber
 Q(H) & := & \int_{R_1 } ^{R_2 } 
\left( 
g^{\prime} (-U^{\prime} ) 
(\Phi ^{\prime} )^2 - f^{\prime } (U) \Phi ^2 \right) 
r^{N-1} \, dr 
\\
\nonumber 
 & = & \int_{R_1 } ^{R_2 }  
g' (-U^{\prime} ) 
(U')^2 (H') ^2 
r^{N-1} \, dr 
+ \int_{R_1 } ^{R_2 } g' (-U' ) (U^{\prime \prime } )^2 
H^2 r^{N-1 }\, dr 
\\
 \nonumber
 & & +
2 \int_{R_1 } ^{R_2 }  g' (-U' ) U^{\prime \prime } U'  H' H r^{N-1 }\, dr 
- \int_{R_1 } ^{R_2 }  f^{\prime } (U) (U')^2 H^2 r^{N-1} \, dr 
\\
\nonumber 
 & =: & Q_1 + Q_2 + Q_3 -Q_4 .
\end{eqnarray}
Since 
$g \in C^1 (0,+\infty )$, 
$f \in C^1 (u_2 , u_1 )$ 
and 
$u(x)= U(|x|)$, 
we have that 
$U\in C^2 (R_1 , R_2 )
\cap C^1 [R_1 , R_2 ]$, 
$U' (R_1 )= U' (R_2 )=0 $ 
and  
\begin{equation}
\label{ODE}
-g' (-U') U^{\prime \prime } r^{N-1 } + (N-1) g(-U' ) r^{N-2} 
= 
f(U) r^{N-1} .
\end{equation}
This yields 
\begin{equation}
\label{calc1}
Q_3  = 
2 (N-1) \int_{R_1 } ^{R_2 } g(-U') U' H' H r^{N-2} \, dr - 
2\int_{R_1 } ^{R_2 } f(U) U' H' H r^{N-1} \, dr.   
\end{equation} 
Further we obtain, using integration by parts, 
\begin{eqnarray*}
Q_2 & = & (N-1) \int_{R_1 } ^{R_2 } g(-U' ) U^{\prime \prime } H^2 r^{N-2} \, dr - \int_{R_1 } ^{R_2 } f(U) U^{\prime \prime } H^2 r^{N-1} \, dr
\\
 & = &
(N-1) \int_{R_1 } ^{R_2 } g (-U') U^{\prime \prime } H^2 r^{N-2} \, dr 
  + \int_{R_1 } ^{R_2 } f '(U) (U')^2 H^2 r^{N-1} \, dr  
\\
 & & + 2\int_{R_1 } ^{R_2 } f(U) U' H' H r^{N-1} \, dr 
+ (N-1 ) \int_{R_1 } ^{R_2 } f(U) U' H^2 r^{N-2} \, dr .
\end{eqnarray*}
Then another integration by parts and (\ref{ODE}) yield
\begin{eqnarray}
\label{calc2}
 & & 
\\
\nonumber
 Q_2 & = & (N-1) \int_{R_1 } ^{R_2 } g(-U') U' H^2 r^{N-3} \, dr  -2 (N-1) 
\int_{R_1 } ^{R_2 }  g(-U') U' H'H  r^{N-2} \, dr
\\
\nonumber
 & & + \int_{R_1 } ^{R_2 } f ' (U) (U')^2 H^2 r^{N-1} \, dr + 2 \int_{R_1 } ^{R_2 } f(U) U' H' H r^{N-1} \, dr. 
\end{eqnarray}
(\ref{QH}) together with (\ref{calc1}), (\ref{calc2}) and assumption (\ref{technass}) give
\begin{eqnarray}
\label{QHnew} & & 
\\
\nonumber
Q(H)  & = & (N-1) \int_{R_1 } ^{R_2 } g(-U' ) U' H^2  r^{N-3 }\, dr + 
\int_{R_1 } ^{R_2 }  g' (-U') (U')^2 (H')^2 r^{N-1} \, dr 
\\
\nonumber 
& \leq & (N-1) \int_{R_1 } ^{R_2 } g(-U' ) U' H^2  r^{N-3 }\, dr + 
\Gamma \int_{R_1 } ^{R_2 }  g (-U') (H')^2 r^{N-1} \, dr 
\\
\nonumber
 & =: & Q_5 (H) + Q_6 (H). 
\end{eqnarray}
An integration of (\ref{ODE}) gives for every $r\in (R_1 , R_2 )$,
\begin{equation} 
\label{intODE}
g(-U' (r) ) r^{N-1}  =  \int_{R_1 } ^r f(U(t)) t^{N-1}\, dt = -\int_{r} ^{R_2 } f(U(t)) t^{N-1} \, dt .
\end{equation}
In view of assumption {\bf (iii)} this means that
\begin{equation}
\label{limQ}
\lim_{r\searrow R_1 } \frac{ g(-U'(r))}{ r-R_1 } = \lim_{r\nearrow R_2 } \frac{ g(-U'(r))}{R_2 -r} =0.
\end{equation} 
Finally, let $\varepsilon \in (0, (R_2 - R_1 )/2 )$, and choose $H=H_{\varepsilon} \in C_0 ^1 (R_1 , R_2 )$, such that $0\leq H_{\varepsilon }  \leq 1 $, $|H_{\varepsilon }' |\leq 2/\varepsilon $ on $(R_1 , R_2 )$ and   
$H_{\varepsilon } (r) \equiv 1 $ for $r\in (R_1 + \varepsilon , R_2 - \varepsilon )$. 
Since $U'<0$, we estimate 
\begin{equation}
\label{Q5}
\lim_{\varepsilon\to 0 } Q_5 (H_{\varepsilon } ) = (N-1 ) \int_{R_1 } ^{R_2 } g(-U' ) U' r^{N-3} \, dr <0 ,
\end{equation}        
and the limits (\ref{limQ}) give 
\begin{eqnarray}
\label{limQ6}
|Q_6 (H_{\varepsilon} )| & \leq & 
\frac{4\Gamma }{\varepsilon ^2 }  \left\{ \int_{R_1 }^{R_1 + \varepsilon} g(-U' ) r^{N-1}\, dr + 
 \int_{R_2-\varepsilon }^{R_2 } g(-U' ) r^{N-1}\, dr 
\right\} 
\\
\nonumber
 & & 
\longrightarrow 0, \quad \mbox{as }\ \varepsilon \to 0.
\end{eqnarray}
Now the assertion follows from (\ref{QHnew}), (\ref{Q5}) and (\ref{limQ6}) by choosing $H= H_{\varepsilon} $ with sufficiently small $\varepsilon >0$.
$\hfill \Box $
\\[0.1cm]
{\bf Remark 4.8.} The assertions of Lemma 4.6. still hold true for $R_1 =0$, that is,  for punctured balls $B_{R_2} \setminus \{ 0\} $. In such situation, the condition $f(u_1 )=0$ in {\bf (iii)}  can be dropped. The proof is analogous and is left to the reader. 
\\[0.1cm]  
{\bf Theorem 4.9.}  {\sl Suppose that $f\in C^1 (0, +\infty )$ and $g$ satisfies assumption {\bf (iv)} of Lemma 4.5. Further, let $u$ be a solution of problem {\bf (P)} which is also a local minimizer  of $J$ in $X$. Then $u$ is radially symmetric and radially decreasing, and moreover, (\ref{rad2}) and (\ref{rad3}) hold. 
}
\\[0.1cm]
{\sl Proof :} $u$ is locally symmetric by Lemma 2.3. Assume that $A_k := B_{R_k } (z_k )\setminus \overline{B_{r_k } (z_k )}$, ($k \in \mathbb{N} $),
 is one of the annuli in the representation (\ref{decomp1}). Then there is a function $U \in C^1 [r_k  , R_k ]$ such that $u(x) = U (|x-z_k |)$ in $A_k $,  $U '(r)<0 $ for $r\in (r_k , R_k )$. Moreover, setting $u_{k} ^- := U(R_k)$   and $u_{k}^+ := U(r_k )$, we have 
$u \geq u_k ^+ $ in $B_{r_k } (z_k )$. It follows that $U\in C^2 (r_k , R_k )$ and $U$ satisfies 
(\ref{ODE}).
\\
Assume first that $R>R_k >r_k \geq 0$. Then we have by Remark 2.2.  that   $U'(R_k )=U'(r_k )=0$, $f(u_k ^- )=0$, and moreover, $f(u_k  ^+ )=0$ in case that $r_k >0$. Since $u$ is a local minimizer to $J$ in $X$, we have that
$J(u+ \varepsilon \varphi ) \geq J(u)$ for all $\varphi \in C_0 ^1 (B_R )$   whenever $|\varepsilon| $ is small enough. This implies
\begin{equation}
\label{2ndvariation} 
\int_{B_R } \left\{ g' (|\nabla u|) \frac{(\nabla u \cdot \nabla \varphi )^2 }{|\nabla u|^2 } + g(|\nabla u|) \left[ \frac{|\nabla \varphi|^2 }{|\nabla u|^2 } - \frac{(\nabla u \cdot \nabla \varphi )^2 }{|\nabla u|^3} \right]   - f ' (u) \varphi ^2 \right\} \, dx \geq 0. 
\end{equation}
In particular, if $\varphi $ has compact support in $A_k $ and is radial, that is,
$\varphi (x) = \Phi (|x-z_k |)$ for some function $\Phi \in C^1 [r_k, R_k]$, (\ref{2ndvariation}) implies
\begin{equation}
\label{2ndvariationrad}
\int_{r_k }^{R_k } \left\{ g' (-U' ) (\Phi ' )^2 - f'(U) \Phi ^2 \right) r^{N-1} \, dr \geq 0.
\end{equation}
But this contradicts to Lemma 4.6.
Hence we must have that $R_k=R$, which means $z_k =0$ and $m=1$. The assertion follows.  
$\hfill \Box $
\\[0.2cm]
{\bf Acknowledgement: } This work was supported by Leverhulme Trust, ref. VP1-2017-004.   

\end{document}